\documentclass[11pt]{amsart}
\usepackage{latexsym}
\usepackage{amssymb,amsmath,amsthm,amscd, 
amssymb,enumerate, verbatim, stmaryrd}
\numberwithin{equation}{section}
\def\co{\colon\thinspace}

\def\G{\Gamma}

\def\d{\delta}
\def\r{\rho}

\def\R{\mathbb{R}}
\def\TF{\mathbb{TF}}

\def\Z{\mathbb{Z}}
\newcommand{\g}{\gamma}

\newtheorem{thm}{Theorem}[section]
\newtheorem{cor}[thm]{Corollary}

\newtheorem{defn}[thm]{Definition}

\newtheorem{Example}[thm]{Example}
\newenvironment{ex}{\begin{Example}\rm}{\end{Example}}
\newtheorem{remark}[thm]{Remark}
\newenvironment{rmk}{\begin{remark}\rm}{\end{remark}}
\newtheorem{Fact}[thm]{Fact}

\newtheorem{Main Lemma}[thm]{Main Lemma}

\newtheorem{Convention}[thm]{Convention}

\begin{document}\abovedisplayskip=6pt plus3pt minus3pt 
\belowdisplayskip=6pt plus3pt minus3pt
\title{\bf
Endomorphisms of relatively hyperbolic groups\rm}
\author[\bf{Igor Belegradek and Andrzej Szczepa{\'n}ski}]
{Igor Belegradek \and Andrzej Szczepa{\'n}ski\\
\\
with appendix by Oleg V. Belegradek}
\date{}
\rm

\address{Igor Belegradek\\School of Mathematics\\ Georgia Institute of
Technology\\ Atlanta, GA 30332-0160}\email{ib@math.gatech.edu}
\address{Andrzej Szczepa{\'n}ski\\ Uniwersytet Gdanski\\ 
Instytut Matematyki\\Wita Stwosza 57\\80-952 Gdansk, Poland}
\email{matas@paula.univ.gda.pl}

\address{Oleg V. Belegradek\\Department of Mathematics\\
Istanbul Bilgi University\\
80370 Dolapdere--Istanbul, Turkey}
\email{\tt olegb@bilgi.edu.tr}

\begin{abstract} 
We generalize some results of Paulin and Rips-Sela 
on endomorphisms of hyperbolic groups to 
relatively hyperbolic groups, and in particular
prove the following.\newline
$\bullet$
If $G$ is a non-elementary
relatively hyperbolic group with slender
parabolic subgroups, and either $G$ is not co-Hopfian
or $\mathrm{Out}(G)$ is infinite, then
$G$ splits over a slender group.\newline
$\bullet$
If $H$ is a non-parabolic subgroup of a  
relatively hyperbolic group, and if any isometric $H$-action 
on an $\R$-tree is trivial, then $H$ is Hopfian.\newline
$\bullet$
If $G$ is a non-elementary relatively hyperbolic group 
whose peripheral subgroups are finitely generated, then $G$
has a non-elementary relatively hyperbolic 
quotient that is Hopfian.
\newline
$\bullet$ 
Any finitely presented group is isomorphic to a finite index 
subgroup of $\mathrm{Out}(H)$ for some group $H$ with Kazhdan property (T).
(This sharpens a result of Ollivier-Wise).
\end{abstract} 
\maketitle

\section{Introduction}
 
The Bestvina-Paulin method~\cite{Pau-gro, Pau-out, Bes-deg}, 
further developed by Sela, has been 
a key ingredient in much of recent work on endomorphisms 
of hyperbolic groups~\cite{Pau-out, RS, Sel-iso, Sel-Kle, Sel-hop},
on endomorphisms of geometrically finite 
Kleinian groups~\cite{Sel-Kle, OP, DP}, and 
in Sela's work on Tarski's problem~\cite{Sel-icm}. 

Very recently, Groves~\cite{Groves} generalized various
results on endomorphisms of hyperbolic groups, notably the Sela's 
shortening argument, to relatively hyperbolic 
groups with finitely generated free abelian parabolic subgroups. 

The point of this paper is to show that some of the easier
applications of the Bestvina-Paulin method extend, with
little effort, to relatively hyperbolic groups
with fairly arbitrary parabolic subgroups. By ``easier'' we mean
the results that do not require JSJ-decompositions or the shortening 
argument.

Relatively hyperbolic groups were introduced by 
Gromov~\cite{Gro-hgr}, and in this paper we use the 
following version of Gromov's definition developed 
by Bowditch (see~\cite[Definition 1]{Bow-rel}).

\begin{defn}\label{defn: rel hyp}\rm
Let $G$ be a group with a \textup{(}possibly empty\textup{)} 
family of subgroups $\mathcal G$.
We say that $G$ is {\it hyperbolic relative to $\mathcal G$} 
if $G$ acts properly discontinuously and isometrically 
on a proper, geodesic, hyperbolic metric space $X$
so that the induced $G$-action on the ideal boundary of $X$ 
is the action of a geometrically finite convergence group 
whose maximal parabolic subgroups are precisely 
the elements of $\mathcal G$, and whose limit set is
the ideal boundary of $X$. Elements of
$\mathcal G$ are called {\it peripheral subgroups}.
\end{defn}

Other definitions of relatively hyperbolic groups 
were developed by Farb~\cite{Far-rel}, 
Bowditch \cite[Definition 2]{Bow-rel}, Yaman~\cite{Yam-rel}, 
Dru{\c{t}}u-Osin-Sapir~\cite{DOS}, Osin \cite{Osi-rel},
Dru{\c{t}}u~\cite{Dru-rel}, and Mineyev-Yaman~\cite{MY}. 
It is known that all these definitions are 
equivalent to Definition~\ref{defn: rel hyp},
provided 
\begin{center}
$G$ and all its peripheral subgroups are finitely 
generated and infinite, 
\end{center}
in which case we say that $G$ {\it satisfies
condition (E).} The proofs of various equivalences can be
found in~\cite[Appendix A]{Dah-rel} (cf.~\cite{Bow-rel, Szc, Bum-rel}),
\cite[Theorem 7.10]{Osi-rel}, \cite{Yam-rel}, \cite[Theorem 8.5]{DOS},
\cite[Theorems 4.21,4.34]{Dru-rel}, \cite[Theorem 57]{MY}.
 
We take this opportunity to correct a slight inaccuracy 
in~\cite[Theorem 7.10]{Osi-rel} and~\cite{Bum-rel} where the
condition (E) is stated without requiring that the peripheral subgroups
are infinite.
In fact, in Definition~\ref{defn: rel hyp} 
each peripheral subgroup is necessarily infinite 
(except when $\mathcal G=\emptyset$) because
parabolic subgroups of a convergence group are infinite.
By contrast, the definitions in~\cite{Far-rel, Osi-rel, DOS, Dru-rel, MY}
allow elements of $\mathcal G$ to be finite, and are
equivalent without requiring that elements of $\mathcal G$ are 
infinite.
It turns out that allowing finite peripheral subgroups does not 
enlarge the class of relatively hyperbolic groups, namely,
Osin proved (for his definition) that if $G$ is hyperbolic 
relative to $\mathcal G$, then for any conjugacy 
invariant subcollection $\mathcal F\subset\mathcal G$ 
of finite subgroups, $G$ is  hyperbolic relative 
to $\mathcal G\setminus\mathcal F$~\cite[Theorem 2.40]{Osi-rel};
thus by making $\mathcal G$ smaller, one can arrange
that either all peripheral subgroups are infinite, or else
$\mathcal G=\emptyset$.

Note that the group $G$ in Definition~\ref{defn: rel hyp}
satisfies (E) if and only if
each peripheral subgroup of $G$ is finitely generated.
Indeed, in this case Definition~\ref{defn: rel hyp} is equivalent 
to~\cite[Definition 2]{Bow-rel}, which means in particular
that $G$ acts simplicially on a connected graph such that
all vertex stabilizers are finitely generated, and 
the quotient graph is finite, which easily implies
that $G$ is finitely generated
(in fact a much more general result 
is proved in~\cite[Theorem 2.2]{Bro}).

Very recently Hruska announced
that Definition~\ref{defn: rel hyp} is equivalent to Osin's definition 
even when $G$ or its peripheral subgroups are not finitely generated. 
This would give many examples of infinitely generated groups $G$
satisfying Definition~\ref{defn: rel hyp}.
However, the only results of this paper that do not need
$G$ to be finitely generated are Theorem~\ref{thm: main}, and 
Corollaries~\ref{cor: fr-hopf},~\ref{cor: fr-endom}.
 
We refer to~\cite[Section 6]{Bow-rel}, or~\cite{Bow-conv}, 
or~\cite[Section 5]{Yam-rel} for relevant background on 
convergence groups. 
A subgroup of a relatively hyperbolic group $G$ is called
{\it elementary} if the limit set of its action as a 
convergence group on the Bowditch's
boundary of $G$ contains at most two points,
which happens exactly if the subgroup is finite, virtually-$\Z$,
or parabolic. Otherwise, the subgroup is called
{\it non-elementary}.
The Tits Alternative for convergence groups~\cite{Tuk}
implies that any small subgroup of a relatively hyperbolic group
is elementary, where a subgroup is called 
{\it small} if it contains no non-cyclic free subgroup.

Examples of relatively hyperbolic groups are:
\begin{itemize}
\item the free products of finitely many finitely generated 
groups are hyperbolic relative to the factors (because 
the free product action on the corresponding
Bass-Serre tree satisfies the Bowditch's 
definition~\cite{Bow-rel} of a relatively
hyperbolic group). 
\item 
hyperbolic groups are hyperbolic relative
to any conjugacy invariant collection $\mathcal G$ of quasi-convex
subgroups such that any element of $\mathcal G$ is equal to its normalizer, and
any two distinct elements of $\mathcal G$ have finite 
intersection~\cite[Theorem 7.11]{Bow-rel}. Any hyperbolic group
is hyperbolic relative to the empty family of subgroups.
\item
geometrically finite isometry groups of Hadamard 
manifolds of negatively pinched sectional curvature
as defined in~\cite{Bow-gf} are hyperbolic relative to
the maximal parabolic subgroups. This includes complete
finite volume manifolds of negatively pinched 
sectional curvature. 
\item 
Osin~\cite{Osi-sc} developed small cancellation theory
for relatively hyperbolic groups, and proved that
small cancellation quotients of relatively hyperbolic groups
are hyperbolic relative to the images of maximal 
parabolic subgroups. These methods imply~\cite{AMO} that 
any two non-elementary finitely generated 
relatively hyperbolic group have a common quotient
that is non-elementary relatively hyperbolic; in particular,
since there are hyperbolic groups with Kazhdan 
property (T), one can arrange the common quotient
to have property (T).
\item
combinations theorems for relatively hyperbolic groups were
proved in~\cite{Dah-comb, Osi-comb}, e.g.
the amalgamation of relatively hyperbolic groups over 
parabolic subgroups 
is relatively hyperbolic, when the parabolic subgroup
is maximal in at least one of the factors.
\item
Hruska-Kleiner~\cite{Hru-Kle} 
proved that $\mathrm{CAT}(0)$-groups with isolated flats
are hyperbolic relative to the flat stabilizers. Examples of
$\mathrm{CAT}(0)$-groups with isolated flats are listed 
in~\cite{Hru-Kle}.
\item Sela's limit group are hyperbolic relative to non-cyclic
maximal abelian subgroups. This was first proved in~\cite{Dah-comb}
using a combination theorem, and in fact according to~\cite{A-Best}
limit groups are $\mathrm{CAT}(0)$-groups with isolated flats.
\end{itemize}

The main observation is as follows.

\begin{thm}\label{thm: main}
Let $H$ be a finitely generated group and $G$ be
a non-elementary relatively hyperbolic group.\newline
$\mathrm{(i)}$
If $\r_k\co H\to G$ is an arbitrary sequence of
pairwise $G$-non-conjugate 
homomorphisms whose images are not parabolic,
then $H$ admits a nontrivial action on an $\R$-tree $T$.\newline
$\mathrm{(ii)}$
Furthermore, if all $\r_k$ are injective, and
$H_a$ is an arbitrary arc stabilizer of the $H$-action on $T$,
then each finitely generated subgroup of $H_a$
is isomorphic to an elementary subgroup of $G$. 
\end{thm}

\section{Proof of Theorem~\ref{thm: main}}

\begin{proof}
Let $(X,d)$ be a proper, geodesic, $\d$-hyperbolic metric space
from the definition of a relative hyperbolic group.
Choose a fundamental domain $F$ for the $G$-action on $X$.
Fix a finite generating
set $S$ for $H$. For $x\in X$ we let 
\[\mu_k(x)=\underset{s\in S}{\mathrm{max}}\ d(x,\r_k(s)(x)),\] and let 
$\mu_k:=\inf_{x\in X} \mu_k(x)$. Choose $x_k\in X$ with
$\mu_k(x_k)\le \mu_k+\frac{1}{k}$. 

First we show that $\mu_k$ has no bounded subsequence.
Arguing by contradiction, consider a bounded subsequence,
which we still denote $\mu_k$. Let $R=1+\sup_{k} \mu_k$.
Choose a $G$-invariant system of horoballs in $X$
such that any two horoballs are at least $R$ apart. 
There exists a sequence $g_k\in G$ such that
$g_k(x_k)\in F$. If the sequence $g_k(x_k)$ is precompact,
then since $G$ acts properly discontinuously,
$g_k\r_k g_k^{-1}$ fall into finitely many 
$G$-conjugacy classes, which contradicts the assumption.
If $g_k(x_k)$ is not precompact, then after passing to 
a subsequence $g_k(x_k)$ lie in a horoball $\mathcal H$,
and since horoballs are at least $R$
apart, $g_k\r_k(s)g_k^{-1}(g_k(x_k))\in \mathcal H$
for each $s\in S$. Hence $g_k\r_k(H)g_k^{-1}$
lies in the parabolic subgroup that stabilizes $\mathcal H$,
which contradicts the assumption that $\r_k(H)$
is not parabolic. 

Then by the work of Bestvina-Paulin 
(see e.g. Sections 3.1-3.5 of~\cite{Bes-surv}), 
the pointed spaces $(\frac{1}{\mu_k} X, x_k)$
subconverges to an $\R$-tree $T$, and $\r_k$'s converge to
a non-trivial isometric $H$-action on $T$, which proves (i). 

Assume $\r_k$ are all injective.
Let $a$ be an arc of $T$ with the stabilizer $H_a$ in $H$.
Let $\check H_a$ be the subgroup of $H_a$ of index at most $2$
that fixes $a$ pointwise. 

Approximate $a$ by geodesic segments $a_k$ in $\frac{1}{\mu_k} X$,
denote the midpoint of $a_k$ by $m_k$.
Changing $\r_k$ within its $G$-conjugacy class,
we can assume $m_k\in F$.

By~\cite[p. 341]{Pau-out} (cf.~\cite[p. 284]{BS}), there exists
a constant $C(\d)$, such that
for each $h,f\in \check H_a$ and all large $k$, 
we have $d(m_k, \r_k([h,f])(m_k))< C(\d)$. 
It follows that for large enough $k$,
the same inequality holds for any given finite 
collection of commutators in $\check H_a$.

First, we prove that if $m_k$ is precompact in $F$, 
then $H_a$ is small, hence elementary.
Indeed, since $m_k$ is precompact,
only finitely elements of $G$ can satisfy 
$d(m_k, g(m_k))< C(\d)$, because
$G$ acts isometrically and properly discontinuously.
We denote the number of such elements by $M$.
If $H_a$ is not small, then $\check H_a$ contains a rank two 
free subgroup generated by $h_1$, $h_2$.
By the previous paragraph, for all large enough $k$,  
$d(m_k, \r_k([h_1,h_2^s])(m_k)< C(\d)$ for $s=1, \dots, M+1$.
So $\r_k([h_1,h_2^{s_1}])=\r_k([h_1,h_2^{s_2}])$ for some $s_1\neq s_2$.
Hence $\r_k([h_1,h_2^{s_1-s_2}])=1$, which is impossible since
$\r_k$ is injective and $h_1,h_2$ generate a free subgroup.

It remains to consider the case when the
sequence $m_k$ is not precompact. 
Fix a $G$-invariant system of horoballs in $X$
such that any two horoballs are at least $C(\d)+1$ apart.
Passing to a subsequence
we can assume that all $m_k$ lie in one horoball.
Denote by $P$ the maximal parabolic subgroup of $G$ stabilizing
this horoball. Since any two horoballs are at least $C(\d)+1$ apart,
and since elements of $G$ map horoballs to horoballs,
the above mentioned result~\cite[p. 341]{Pau-out} implies that
the $\r_k$-image of any finite collection of commutators lies 
in $P$, provided $k$ is large enough. 

Let $z\in\partial X$ 
be the unique fixed point of $P$. Note that $P$ is the stabilizer
of $z$ in $G$.
Let $\G$ be an arbitrary finitely generated subgroup of $H_a$,
and let $\check\G:=\G\cap \check H_a$; the
index of $\check\G$ in $\G$ is $\le 2$, so
$\check\G$ is a normal subgroup of $\G$.
If all commutators in $\check\G$ have finite order (i.e.
the commutator subgroup $[\check\G, \check\G]$ is torsion), 
then $\check\G$ and $\G$ must be small, 
hence $\r_k(\G)$ is elementary.

Suppose that there exists a commutator $q\in\check\G$ of
infinite order. Fix an arbitrary $s\in\check \G$, and 
take $k$ sufficiently large so that the $\r_k$-image of the commutators 
$sqs^{-1}q^{-1}$ and $q$ lies in $P$, hence 
$\r_k(sqs^{-1})$ lies in the intersection of
$P$ and $\r_k(s)P\r_k(s^{-1})$. 
Since each $\rho_k$ is injective,
$\rho_k(q)$ has infinite order, hence the intersection of   
$P$ and $\r_k(s)P\r_k(s^{-1})$ is infinite. 
Distinct maximal parabolic subgroups have finite
intersection, so $P=\r_k(s)P\r_k(s^{-1})$, i.e. $\r_k(s)\in P$. 
Since $\check\G$ is finitely generated, 
$\r_k(\check\G)\subset P$ if $k$ is large enough.
To see that $\r_k(\G)$ lies in $P$
take an arbitrary $\g\in\r_k(\G)$ and note that 
the subgroup $\g\r_k(\check\G)\g^{-1}$ fixes both
$\g(z)$ and $z$, because
$\check\G\trianglelefteq\G$. Since 
the group $\g\r_k(\check\G)\g^{-1}$ is infinite parabolic, 
it has a unique fixed point, so $\g(z)=z$. Then
$\r_k(\G)$  lies in the stabilizer of $z$,
which equals to $P$.
\end{proof}

\begin{rmk} It follows from Theorem~\ref{thm: main}(ii) that
if parabolic subgroups of $G$ are
small, then each $H_a$ is small. The same is true
if small is replaced by amenable, slender, finitely generated
virtually nilpotent, finitely generated virtually abelian, etc. 
Also if parabolic subgroups of $G$ are free, and
$G$ is torsion-free, then each $H_a$ is locally free.
\end{rmk}

\section{Endomorphisms and actions on $\R$-trees}

Theorem~\ref{thm: main} immediately implies the following 
corollary, which for hyperbolic groups 
is due to Paulin~\cite{Pau-out}.

\begin{cor}\label{cor: out}
If a finitely generated non-elementary 
relatively hyperbolic group $G$ has 
infinite $\mathrm{Out}(G)$, then
$G$ acts non-trivially on an $\R$-tree such that
any finitely generated subgroup of each arc stabilizer 
is isomorphic to an elementary subgroup of $G$.
\end{cor}

A group is called {\it co-Hopfian}
if every injective endomorphism of the group is surjective.

\begin{cor} \label{cor: cohopf}
Suppose $G$ is a finitely generated non-elementary relatively 
hyperbolic group that contains no infinite torsion group
and is not isomorphic to a parabolic subgroup of $G$.
If $G$ is not co-Hopfian, then 
$G$ acts non-trivially on an $\R$-tree such that
any finitely generated subgroup of each arc stabilizer 
is isomorphic to an elementary subgroup of $G$.
\end{cor}
\begin{proof} Arguing by contradiction, we wish to show
that $G$ is co-Hopfian.
In the proof of~\cite[Theorems 3.1]{RS} 
Rips-Sela establish the co-Hopf property for any group $G$ 
such that \newline
(a) $G$ has no infinite torsion subgroup,\newline
(b) up to conjugation there are only finitely many
monomorphisms $G\to G$,\newline
(c) the image of any monomorphism $r\co G\to G$ has finite
centralizer in $G$.

Now (a) is true by assumption, (b) follows from 
Theorem~\ref{thm: main} and the contradiction assumption.
To verify (c) note that $r(G)$ is
non-elementary because we assumed
G is non-elementary, so in particular, it is not finite or
virtually-$\Z$,
and also by assumption G is not isomorphic to a parabolic subgroup.
As we note in the proof of Corollary~\ref{cor: fr-endom},
any non-elementary subgroup of a relatively hyperbolic group
has finite centralizer.
\end{proof}

\begin{rmk}
Sela~\cite{Sel-Kle} proved that
a torsion-free hyperbolic group is co-Hopfian if and only if
it is freely indecomposable. Delzant-Potyagailo characterized
co-Hopfian geometrically finite Kleinian groups in terms of 
splittings~\cite{DP}.
\end{rmk} 

\begin{rmk}
We do not know whether the assumption in Corollary~\ref{cor: cohopf}
``$G$ contains no infinite torsion group'' can be dropped.
What is actually used in the proof is that if $\phi$ is a
non-surjective monomorphism of $G\to G$ and $A_{k,\phi}$ is the
(necessarily finite) centralizer of the subgroup $\phi^k(G)$, 
then the torsion subgroup $A_\phi=\bigcup_k A_{k,\phi}$ is finite.
\footnote{It was recently 
showed in~\cite[Lemma 4.44]{DS} that $A_\phi$ is finite if $G$ 
satisfies (E).}
The assumption ``$G$ is not isomorphic to a parabolic subgroup of $G$''
in Corollary~\ref{cor: cohopf} holds e.g. if parabolic subgroups are 
small and $G$ is non-elementary,
or if parabolic subgroups are co-Hopfian.
This assumption cannot be dropped as is shown in the following example
(recall that groups with property (T) have no nontrivial actions
on $\R$-trees).
\end{rmk}

\begin{ex}\it \label{ex: not-cohopf-T}
There exists a non-elementary torsion-free
relatively hyperbolic group that is not co-Hopfian, 
and that has Kazhdan property (T).
\end{ex}
\begin{proof}
Osin's small cancellation methods allow to construct
relatively hyperbolic Kazhdan groups with prescribed
maximal parabolic subgroups.
Specifically, if $H$ be an infinite finitely generated group 
and $K$ is a torsion-free non-elementary hyperbolic Kazhdan group, then
$G:=H\ast\Z\ast K$ is hyperbolic relative to the factors
$H$, $\Z$. It follows from properties of free products
that $K$ is a suitable subgroup of $G$ (see~\cite{Osi-sc}
for a definition). By~\cite[Theorem 2.4]{Osi-sc}
there is an epimorphism $\eta\co G\to \bar G$ such that 
$\eta(K)=\bar G$, the restriction of $\eta$ to $H\cup \Z$
is injective, $\bar G$ is (Osin) hyperbolic relative to 
$\eta(H)$, $\eta(\Z)$,
and $\bar G$ is obtained from $G$ by adding finitely many 
relations. In particular, $\bar G$ is Kazhdan (because property (T) is
inherited by quotients), $\bar G$ is non-elementary,
and furthermore if $H$ is finitely presented, 
then so is $\bar G$. Moreover, by~\cite[Theorem 2.4]{Osi-sc}
any element of finite order in $\bar G$ is the image of
a finite order element in $G$, so if $H$ is torsion-free,
then so are $G$ and $\bar G$. By construction, $G$ satisfies (E)
so it is relatively hyperbolic in the sense of 
Definition~\ref{defn: rel hyp}. 
Now suppose $H$ is a universal torsion-free
finitely presented group (see Theorem~\ref{thm: univ tf} below), 
i.e. $H$ is a finitely presented 
torsion-free group such that any finitely presented torsion-free 
group embeds into $H$.   
Thus $\bar G$ embeds into $H\cong\eta(H)\subset\bar G$, hence
$\bar G$ is not co-Hopfian.
\end{proof}

\begin{rmk}  
Corollaries~\ref{cor: out},~\ref{cor: cohopf} generally
fail for
finitely generated subgroups of hyperbolic groups:
a non-co-Hopfian finitely generated subgroup of a 
hyperbolic group was constructed in~\cite{KW},
and by~\cite{OW} there exist Kazhdan groups that have 
infinite outer automorphism groups yet are embedded 
into hyperbolic group. However, the examples 
in~\cite{KW, OW} are not finitely presented, and
it seems no finitely presented examples are known.
\end{rmk}

\section{Splitting over slender groups}

A group is called {\it slender} 
if all its subgroups are finitely generated.
The class of slender group is closed under
extensions~\cite{Dun-S}, in particular, virtually polycyclic 
groups are slender. Slender groups are small
because non-abelian free groups contain infinitely 
generated subgroups. 

As explained in the introduction, 
if $G$ is a  relatively hyperbolic group 
with slender (or more generally finitely generated)
peripheral subgroups, then $G$ is finitely generated.
Furthermore, by the Tits Alternative for relatively 
hyperbolic groups~\cite{Tuk},
small subgroups of $G$ are parabolic, hence slender. 
Then since in any group an ascending
chain of small subgroups is small, every 
ascending chain of slender subgroups of $G$ stabilizes.
Hence, any nontrivial
minimal $G$-action on an $\R$-tree with slender
arc stabilizers is stable~\cite{BF}, and therefore~\cite{BF} if
$G$ is finitely presented, then there exists a 
nontrivial $G$-action on a simplicial tree slender edge 
stabilizers, i.e. a splitting of $G$ over a slender subgroup. 
(Note that the splitting is nontrivial provided $G$ is 
non-elementary). 
In fact, according to Dunwoody~\cite[Theorem 2]{Dun}, 
the assumption ``$G$ is finitely presented'' can be replaced
by ``$G$ is finitely generated'', which holds if $G$ is relatively
hyperbolic with slender peripheral subgroups.
Thus combining~\cite[Theorem 2]{Dun} with 
Corollaries~\ref{cor: out},~\ref{cor: cohopf}, 
we get the results stated in the abstract about co-Hopf property
and infinite outer automorphism group, as well as
the following.

\begin{cor}\label{cor: slender}
Let $H$ be a finitely generated group that is not slender, 
and let $G$ be a relatively hyperbolic group with slender 
parabolic subgroups.
If $\r_k\co H\to G$ is an arbitrary sequence of
pairwise $G$-non-conjugate injective
homomorphisms, then $H$ admits a nontrivial 
splitting over a slender group.
\end{cor}

\begin{rmk} It is interesting to compare
Corollary~\ref{cor: slender} 
with a (much more delicate) 
result of Dahmani~\cite{Dah-acc} who showed
that if $H$ is a finitely presented group and 
$G$ is a relatively hyperbolic group, then up to $G$-conjugacy
there are only finitely many subgroups of $G$ that are
non-parabolic, do not split over parabolic subgroups, and are
the images of homomorphisms $H\to G$. 
\end{rmk}

\section{Subgroups with property $F\R$ are Hopfian}

A group $H$ is said to have
{\it property $F\R$} if any 
$H$-action on an $\R$-tree fixes a point, i.e. is not 
{\it non-trivial}. Any group with property $F\R$ is finitely 
generated~\cite[p.81]{Ser}.
By Theorem~\ref{thm: main} if $H$ 
has property $F\R$, then there are only finitely many 
$G$-conjugacy classes
of homomorphisms $H\to G$ with non-parabolic images.

The class of groups with property 
$F\R$ is closed under extensions~\cite[p. 85]{Ser} and quotients, 
and contains 
\begin{itemize}
\item groups with Kazhdan property (T) (this was proved in~\cite{Wat} for 
simplicial trees, and extended to $\R$-trees in~\cite{Nos}); 
\item $\mathrm{Aut}(F_n)$ and $\mathrm{Out}(F_n)$ for $n>2$~\cite{CV},
where $F_n$ is a free group of rank $n$; 
mapping class groups $M_{g,r}$ for $g\ge 2$~\cite{CV};
\item 
finitely generated 
Coxeter groups with no $\infty$ labels in their 
Coxeter diagrams (this is an exercise in~\cite[p. 93]{Ser}), 
\item fundamental groups of closed irreducible non-Haken 
$3$-manifolds~\cite[Proposition 2.1]{Mor-Sha}.
%
%
\end{itemize}

Recall that group is called {\it Hopfian}
if every surjective endomorphism of the group is injective.

\begin{cor}\label{cor: fr-hopf}
If a non-parabolic subgroup $H$
of a relatively hyperbolic group
has property $F\R$, then $H$ is Hopfian. 
\end{cor}
\begin{proof}
If $\r$ is a non-bijective epimorphism of $H$, then 
precomposing the inclusion $H\hookrightarrow G$ with 
the powers $\r^k$, we get 
the homomorphisms $H\to G$
with image $H$, which is non-parabolic by assumption.
To show that $\r$ is injective, it is enough to prove 
that some power of $\r$ is injective. Look at the sequence 
$\psi_k:=\r^{2^k}$, so that
for any two endomorphisms in the sequence $\{\psi_k\}$
the one with the larger index is a power of the other.
By Theorem~\ref{thm: main}, after passing to a subsequence,
we can assume that $\psi_k$'s are conjugate in $G$. 
In particular, if $\psi$ denotes the first endomorphism
in the subsequence, we get $\psi^{s_k}=i_{g_k}\circ \psi$, 
where $g_k\in G$ and $i_{g_k}$ is the corresponding 
inner automorphism of $G$.
Then $\psi^{s_k}$ and $\psi$ have equal kernels.
Take $\g\in \ker(\psi^{s_k-1})$.
Since $\psi$ is onto, we can find $\tilde\g\in H$
with $\psi(\tilde\g)=\g$. So 
$\tilde \g\in\ker(\psi^{s_k})=\ker(\psi)$, which implies $\g=1$, 
as wanted.
\end{proof}

\begin{cor}\label{cor: quot-hopf}
If $G$ is non-elementary relatively hyperbolic group 
whose peripheral subgroups are finitely generated, then $G$
has a non-elementary relatively hyperbolic 
quotient $\bar G$ that is Hopfian.
Also if $G$ is hyperbolic, then so is $\bar G$.
\end{cor}
\begin{proof} 
Since the peripheral subgroups are finitely generated, 
$G$ satisfies $(E)$, so Definition~\ref{defn: rel hyp}
is equivalent to Osin's definition for which
it is known by~\cite{AMO} that any two finitely generated
non-elementary relatively hyperbolic groups 
have a common non-elementary relatively hyperbolic quotient.
Since there exists a non-elementary
hyperbolic group with Kazhdan property (T),
any non-elementary relatively hyperbolic group $G$ has a 
non-elementary (Osin) relatively hyperbolic quotient
$\bar G$ with property (T), which is Hopfian by Corollary~\ref{cor: fr-hopf}.
By~\cite[Theorem 2.4]{Osi-sc} (on which the results~\cite{AMO} are based)
the surjection $G\to\bar G$ maps the peripheral 
subgroups of $G$ isomorphically onto the peripheral 
subgroups of $\bar G$, so $\bar G$ satisfies $(E)$
and hence $\bar G$ is also relatively hyperbolic in the sense of
Definition~\ref{defn: rel hyp}. Finally, if $G$ is hyperbolic, we can take 
$\mathcal G$ to be empty, so that $\bar G$ has no peripheral subgroups,
i.e. $\bar G$ is hyperbolic.
\end{proof}

\begin{rmk}
That hyperbolic groups with property $F\R$ are Hopfian 
was noted to~\cite[Theorem 2.1]{RS}, and 
Corollary~\ref{cor: fr-hopf} extends this result to
relatively hyperbolic case. Later 
Sela~\cite{Sel-hop} showed that any 
torsion-free hyperbolic group is Hopfian, and more recently
Bumagin proved~\cite{Bum-hop} that any finitely generated
subgroup of a torsion-free hyperbolic group is Hopfian.
Thus the assumption ``$H$ has property $F\R$'' in 
Corollary~\ref{cor: fr-hopf} is probably far from optimal. 
On the other hand, there certainly exist non-Hopfian 
relatively hyperbolic groups, e.g. the free product of any
two finitely generated groups one of which is non-Hopfian.
\end{rmk}

\begin{rmk}
Ollivier-Wise~\cite{OW} and de Cornulier~\cite{dC} gave examples of
non-Hopfian Kazhdan groups. By Corollary~\ref{cor: fr-hopf}
these groups cannot be embedded
in a relatively hyperbolic group as non-parabolic subgroups, because
Kazhdan groups are finitely generated and have property $F\R$.
\end{rmk}

\section{Automorphism groups of Kazhdan groups}

Any hyperbolic group with Kazhdan property (T) has finite
outer automorphism group~\cite{Pau-out}. By contrast,
Ollivier-Wise~\cite{OW} showed that any finitely presented group
is isomorphic to the quotient of a torsion-free hyperbolic 
group $G$ by an infinite normal subgroup $H$ that has Kazhdan
property (T). Then the canonical homomorphism
$G/H\to\mathrm{Out}(H)$ is injective, thus Ollivier-Wise
concluded that any finitely presented group embeds into 
$\mathrm{Out}(H)$ for some group $H$ with property (T). 
We note that the embedding $G/H\to\mathrm{Out}(H)$ 
necessarily has finite cokernel.

\begin{cor} \label{cor: fr-endom}
If $H$ has property $F\R$ and is a normal subgroup
of a non-elementary torsion-free relatively hyperbolic group 
$G$, then the canonical homomorphisms
$G/H\cong\mathrm{Inn}(G)/\mathrm{Inn}(H)\to \mathrm{Out}(H)$ 
and $G\cong\mathrm{Inn}(G)\to\mathrm{Aut}(H)$ are 
embeddings onto finite index subgroups.
\end{cor}
\begin{proof} 
We think of $G$ as a discrete convergence group.
Infinite elementary subgroups have elementary normalizers
because the normalizer stabilizes the limit set of the subgroup.
Since $H$ is normal and $G$ is non-elementary, 
we conclude that $H$ is non-elementary. 

Now we show that 
non-elementary subgroups have finite centralizers.
Indeed, if the centralizer $C_G(H)$ of $H$ in $G$
contains an infinite order 
element $c$, then $H$ fixes the limit set of the cyclic
subgroup generated by $c$, so $H$ cannot be non-elementary.
Thus $C_G(H)$ must be a torsion group, in particular,
$C_G(H)$ is elementary, else $C_G(H)$ would have to
contain a non-abelian free subgroup. If $C_G(H)$ is infinite,
it has a nonempty limit set that must be fixed by $H$,
contradicting the assumption that $H$ is non-elementary.

Since $G$ is torsion free, $C_G(H)$ is trivial.
Hence the canonical
homomorphism $\mathrm{Inn}(G)\to \mathrm{Aut}(H)$
in injective, and we identify its image with
$\mathrm{Inn}(G)$. 
By Theorem~\ref{thm: main} applied to the inclusion 
$H\hookrightarrow G$ precomposed with automorphisms of $H$,
$\mathrm{Inn}(G)$ has finite cokernel in $\mathrm{Aut}(H)$.
Let $G_0$ be the intersection of all the
conjugates of $\mathrm{Inn}(G)$ in $\mathrm{Aut}(H)$. Note
that $G_0$ is a finite index normal subgroup of $\mathrm{Aut}(H)$
that contains $H$.
The kernel of the surjection 
of $\mathrm{Out}(H)=\mathrm{Aut}(H)/\mathrm{Inn}(H)$ 
onto the finite group $\mathrm{Aut}(H)/G_0$
equals to $G_0/\mathrm{Inn}(H)\cong G_0/H$.
Since $G_0/H\le G/H\le\mathrm{Out}(H)$, 
we conclude that $G/H$ has finite index in $\mathrm{Out}(H)$.
\end{proof}

\begin{cor}\label{cor: kazh-out} Given a finitely presented
group $Q$ there exists a group $H$ with Kazhdan property (T)
such that $Q$ is isomorphic to a finite index
subgroup of $\mathrm{Out}(H)$, and the group $\mathrm{Aut}(H)$ is hyperbolic.
\end{cor}
\begin{proof} By~\cite{OW}, $Q\cong G/H$ where $G$ is hyperbolic and $H$
has property (T). By Corollary~\ref{cor: fr-endom}, $Q$ embeds into
$\mathrm{Out} (H)$ as a finite index subgroup, and
$\mathrm{Aut}(H)$ is hyperbolic because it
contains the hyperbolic subgroup $G$ of finite index.
\end{proof}

\appendix
\section{On universal torsion-free finitely presented groups,  
by Oleg V. Belegradek}
\label{sec: tf-higman}

\begin{thm}\label{thm: univ tf}
There exists a
universal torsion-free finitely presented group.
\end{thm}
\begin{proof}
We will show below that there exists a torsion-free
recursively presented group $P$ which contains an 
isomorphic copy of every torsion-free finitely presented group.
(Note that this is not completely trivial:
the naive idea ``consider an effective listing of
finite presentations of all torsion-free finitely 
presented groups and form the disjoint union of the 
presentations'' fails because such a listing does 
not exist~\cite{Lem}.)
By the Higman Embedding Theorem, $P$ embeds into
a finitely presented group $H$ which is built from $P$
and the trivial group by a finite sequence of
HNN-extensions and free 
products~\cite[pp. 364--365, 389]{Rot-3rd-edition}. 
Since the class $\TF$
of all torsion-free groups
is closed under
HNN-extensions and free products, $H$ is torsion-free.
Thus, $H$ is a universal torsion-free finitely presented group.

The class $\TF$ is a quasivariety
specified by the recursive set $Q$ of quasi-identities
$\forall x\,(x^n=1\to x=1)$, where $n$ runs over positive integers.
By general nonsense~\cite[Theorem 9.2.2]{Hod}, any quasivariety
admits presentations. We denote by $G^\tau$ the group
defined in $\TF$  by a presentation $\tau$,
and by $G_\tau$ the group defined by $\tau$ in the variety
of all groups.
Let $(\pi_n: n=0,1,\dots)$ be an effective listing
of all finite presentations. Clearly, the disjoint union $\pi$
of all $\pi_n$
is a recursively enumerable presentation.
We show that one can take $G^\pi$ as $P$.

Any torsion-free finitely presented group
embeds into $G^\pi$. Indeed, if $\tau$ is
a finite presentation such that $G_\tau$ is torsion-free
then $G_\tau=G^\tau$; clearly, $G^\pi$ contains an 
isomorphic copy of $G^\tau$.
It remains to show that $G^\pi$ is recursively presented
in the variety of all groups. Let $\pi=\langle X,R\,\rangle$,
and $R^\prime$ be the set of all relations of
$G^\pi$ in generators $X$.
Then $w$ is in $R^\prime$ if and only if the formula
$w=1$ is a consequence of the set of formulas
$Q_R=Q\cup\{r=1: r\in R\}$
in first order logic~\cite[Lemma 9.2.1]{Hod}.
Since $Q$ and $R$ are recursively enumerable,
$Q_R$ is recursively enumerable, too.
The set of consequences of any recursively enumerable set
of axioms is recursively enumerable~\cite[Lemma 6.1.3]{Hod}, 
hence $R^\prime$ is recursively enumerable, as claimed.
\end{proof}

\begin{rmk} The above argument applies verbatim to show that
every quasivariety of groups $K$ defined by a recursively enumerable
set of quasi-identities contains a group which is recursively presented in
the variety of all groups and contains isomorphic copies of 
all finitely presented groups
that belong to $K$. For the quasivariety $\TF$ there is
a more constructive (but less elegant)
proof that $G^\pi$ is recursively presented based on the following
observation: $R^\prime=\bigcup_n R_n$, where
$R_0$ is the normal closure
of $R$ (in the free group generated by $X$), and
$R_{n+1}$ is the normal closure
of the set of all roots of all elements of $R_n$.
\end{rmk}

\bf{Acknowledgements:}\rm\
We are grateful to Oleg V. Belegradek,
Fran{\c{c}}ois Dahmani, Chris Hruska, Ilya Kapovich, 
Misha Kapovich, Denis Osin, Mark Sapir,
and the referee for helpful communications.
The first author was partially supported by the NSF grant \# DMS-0352576.
The second author was partially supported by the Polish grant
BW 5100-5-0096-5.

\bf Added in April 2007:\rm\
A version of this paper was posted as~\cite{BelS} 
in January of 2005, and the only new results in the present 
version are Corollary~\ref{cor: quot-hopf} and 
Appendix~\ref{sec: tf-higman} (whose sole purpose was to 
make the group in Example~\ref{ex: not-cohopf-T} torsion-free).
Eleven months later Dru{\c{t}}u-Sapir circulated a 
remarkable preprint~\cite{DS} in which they extended 
Rips theory to tree-graded spaces and, 
in particular, obtained generalizations of Theorem~\ref{thm: main}
and Corollaries~\ref{cor: out},~\ref{cor: cohopf},~\ref{cor: slender}.
Dru{\c{t}}u-Sapir's argument is different and quite long. 

\small
\bibliographystyle{amsalpha}
\bibliography{rh-revised}
\end{document}